# SYSTEMS ENGINEERING OF OPTIMAL CONTROL SYNTHESIS OF THE STRUCTURE OF THE TECHNOLOGICAL PRODUCTS CONVERSION SYSTEM (PART 1)

I. Lutsenko
PhD, Professor
Department of Electronic Devices
Kremenchuk Mykhailo Ostrohradshyi
National University
Pervomaiskaya str., 20,
Kremenchuk, Ukraine, 39600
E-mail: delo-do@i.ua

*Розроблено архітектуру системи управління з можливістю повномасштабної параметричної оптимізації. Встановлено, що це можливо, якщо функції перетворення і буферизації технологічного продукту виконують різні системи. Синтезовано архітектуру системи перетворення на прикладі технологічного процесу нагрівання рідини. Створювані моделі тестувалися і досліджувалися в спеціально розробленому для цих цілей безкоштовному програмному конструкторі EFFLI*

*Ключові слова: синтез систем, керована система, система перетворення, система буферизації, структура системи*

*Разработана архитектура управляемой системы с возможностью полномасштабной параметрической оптимизации. Установлено, что это возможно, если функции преобразования и буферизации технологического продукта выполняют разные системы. Синтезирована архитектура системы преобразования на примере технологического процесса нагрева жидкости. Создаваемые модели тестировались и исследовались в специально разработанном для этих целей бесплатном программном конструкторе EFFLI*

*Ключевые слова: синтез систем, управляемая система, система преобразования, система буферизации, структура системы*

## 1. Introduction

An enormous amount of publications deal with optimization or optimal control problems, and the number of publications continues to grow like an avalanche. Herewith, the fact that the optimal control technology can be fully implemented only if there is a well-defined system capability is overlooked.

For example, stocks of raw and energy products are formed at the input and finished products at the output of almost every enterprise. On the one hand, warehouses, tanks, various repositories, etc. (in general, buffering mechanisms) allow to adjust the volumes of the input batches. Thus, it becomes possible to change transport charges, investment volumes, regulate the safety stock level, and so on.

At the same time, due to the presence of the buffering system, a consumer or the internal consumption system can get instant access to the desired product in the required volume, despite significant fluctuations in demand.

The situation is different for intersystem processes. Lack of the efficiency criterion and high cost of process equipment has led to the established practice of minimizing the equipment outlay and optimization capabilities of controlled systems respectively. As a consequence, fairly sophisticated automated complexes and controlled systems are produced today, most of which are either not adapted to solving optimization problems, or those capabilities are very limited, as a rule, by the process parameters stabilization problems.

At the same time, it is obvious that any system operation is performed with the single purpose of increasing the value of the finished product relative to the value of initial (raw, energy and information) products. This means that the main aspect the developers of controlled systems should pay attention to – is not just to ensure the manufacture of high-quality products, but also to have the ability to change the value added as a result of performing the system operation. Only in this case, full parametric optimization by the resource efficiency criterion can be implemented. This problem especially concerns energy- and resource-intensive controlled systems.

## 2. Analysis of the literature data and problem statement

The basic descriptive document for developers of controlled systems is the MRP, ERP and CRM control standards manual [1]. Analysis of the contents of these documents shows that they do not have a description of the optimal controlled system construction concept and definition of the efficiency criterion.

One of the general trends of systems engineering is the perception of the system as an independent object, the processes of which are subject to optimization [2].

Motor transport, which can operate autonomously for some time can be an example [3]. However, movement time of the vehicle, its trajectory or inner processes can not be subject to optimization. The subject of optimization can be cargo movement operation, where transport acts as a movement mechanism. The driver, the vehicle and the road (overpass, pipeline) is a movement system. Processes of the system can not be optimized as it is often tried to do [4]. It is possible to optimize, for example, the movement process of the system. But in this case, the system is the movement product.

Certainly, having internal stocks of the energy product or access to these stocks, it is possible to control their feed and, thus, change expenditure through the power consumption minimization [5]. But is the power consumption minimization problem the control objective? If so, the system is better not to run. We will obtain a zero consumption of the energy product.

A common optimization objective is process time minimization [6] or selection of the most successful motion path [7]. In this case, the control objective remains aloof again.

The global control objective is to ensure profit maximization within a certain time interval. And what will the manager do if the system operates in a mode that does not correspond to the minimum costs, but increases them for example by 5 % and raises annual profit by 20 % with respect to the austerity measures? Will he reduce the energy consumption by 5 % and annual profit by 20 %?

Practitioners, having the opportunity to find the way to higher incomes are unlikely to specifically go into the theory of optimization by the minimum costs.

This means that the subject of the research is not inside the system under investigation, but outside it. This also means that before solving the parametric optimization problems, it is necessary to solve the controlled system synthesis problem, basically allowing optimization. Otherwise, an attempt of local optimization within a certain system can lead to a global efficiency reduction within the controlled system as a whole.

Similar claims can be made against the control, which is focused on the minimum time of the system operation [8]. It is well known that the performance improvement leads to increased wear of expensive equipment.

Hence, the simple conclusion that the architecture of the controlled system must take its own structure as one of the input technological products.

Practice shows that since "experienced optimizers" can not explain to the industrialist the way to find a compromise between energy saving, equipment wear, operation time and profit, technologist, even if there is a choice of control, simply sets the mode that in his *personal opinion* is the best.

In general, the analysis of publications shows that the authors do not consider the system under investigation as one of numerous objects of the controlled system, between which control signals and technological products, including system equipment as one of the technological products circulate. Furthermore, models of controlled systems do not use comparable cost values as the scaling factors.

Exceptions are unsuccessful attempts to select optimal managerial decisions by the method of direct estimation of value added within a certain time interval based on economic indicators. As the authors note, "this is caused by the complexity of determining the economic indicators with a sufficient degree of accuracy and reliability in short periods of the operation time of the enterprise (shift, day, week)" [9].

Indeed, the direct estimation method is inconvenient for use in real systems because of the complexity of the reproduction of one-type operations, a significant impact of external factors, and, primarily, because of the continuance of the estimation process. It is irreplaceable for determining the adequacy of efficiency criteria in modeling problems [10] and is one of the main tools in evaluating the achievements in the optimal system structure synthesis problems.

## 3. Goal and objectives

The goal of the paper is the development of a conceptually unified architecture of basic optimal control systems that allows achieving the strategic goal of control – the most rapid increase in the value added of the system at a certain time interval due to full parametric optimization taking into account the cost parameters of system products and process equipment wear.

Tasks, the solution of which is necessary to achieve the goal are:
– the synthesis of the structure of the conversion system with batch feed of raw products;
– the development of the controlled system of extreme control with batch feed of raw products;
– the development of systems of direct estimation of the effectiveness of the technological process; synthesis of the separation system structure;
– the optimization of the separation system processes;
– the development of the software constructor for the verification of the performance of systems, mechanisms, structures and objects, as well as for the research of operations and processes of controlled systems.

## 4. Synthesis of the conversion system with batch feed of raw products

The first step in the synthesis of a structure that provides for optimal control, is the development of the basic architecture of the controlled system, which ensures control of internal processes with the required number of degrees of freedom.

The experiments have shown that this freedom can be obtained only if the functions on achieving a given quality and mode of release of the required amount of output product are divided among highly specialized systems. Each such system performs only one system operations. In this case, the function of batch heating is performed by the conversion system, and function of the release of the required volume of the finished product per unit of time – by the buffering system.

Systems with a *batch feed of raw materials* in the scientific and academic literature are commonly called *periodic systems*. However, systems with batch feed of products must be named more correctly since these systems are in the search mode for most of the functioning time and the operation time in this case, and hence the period varies continuously.

As an example of the construction of the controlled system with the necessary degrees of freedom, let us consider the controlled liquid heating system. This choice is caused by the heating system inertia, heating model simplicity and the possibility of the analytical determination of the electric heating mechanism wear depending on the control mode.

The controlled liquid heating system consists of the cold liquid feed system, energy product feed system, heating system and heated liquid buffering system (Fig. 1).

Following the feed of the setpoint signal $z_{wp}$, associated with the need of the heated liquid buffering system restocking, heating system generates a control signal $u_w$ to feed a certain amount of cold liquid $r_w$. Completion of the cold liquid feed process leads to the formation of a control signal $u_p$, which provides the feed of the energy product $r_p$ with a predetermined intensity. This naturally leads to an increase in the heated liquid temperature. Once the liquid temperature reaches a predetermined value $z_T$, a control signal $u_{off}$ to shut down the energy product feed is generated, and the heated liquid $u_{off}$ is transferred to the buffering system.

Providing a consumer with finished products is carried out by the consumption signal $z'_{wp}$. An appropriate amount of the finished product per unit of time $p'_w$ is transferred to the consumer from the buffering system.

Heating system operates in a cyclic mode, and can not independently provide the consumer with a continuous flow of the heated liquid. However, it allows to independently select the liquid portion heating intensity that allows to control the heating rate (energy product expenditure, heating mechanism wear), and thus choose the most convenient heating mode based on the higher grounds.

In turn, the heated liquid buffering system provides consumer with a finished product that satisfies him by both qualitative and quantitative indicators.

Furthermore, the buffering system has the ability of independent regulation of both safety stock level and upper level of its buffering mechanism.

For the synthesis of the internal structure of the conversion system, let us map the interaction of its objects, which provide data processing and information exchange, as a Petri net (Fig. 2). Simple transitions are not mapped in order to simplify the perception. Each element of the CS represents one simple operation.

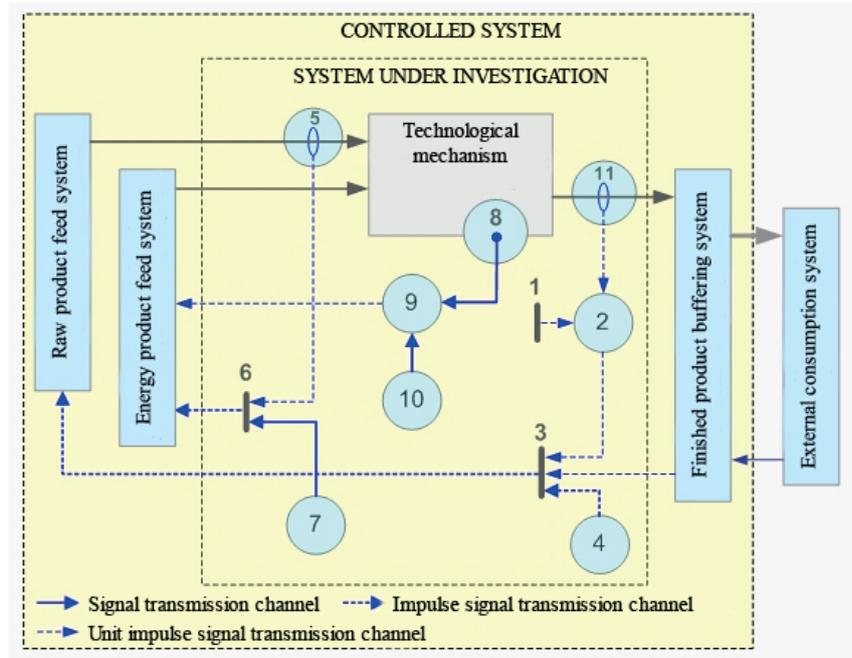

Fig. 2. Controlled liquid heating system with detailing of the synthesized heating system: 1 — system startup; 2 — multiplexer; 3 — data feed synchronization mechanism; 4 — object of formation-feed of setting signal relative to the cold liquid volume; 5 — registration mechanism of the time of the cold liquid feed completion; 6 — data feed synchronization mechanism; 7 — object of formation-feed of setting signal relative to the energy product intensity; 8 — determination mechanism of the current value of the heated liquid temperature; 9 — comparator; 10 — feed mechanism of reference value of the heated liquid temperature; 11 — registration mechanism of the time of the heated liquid feed completion

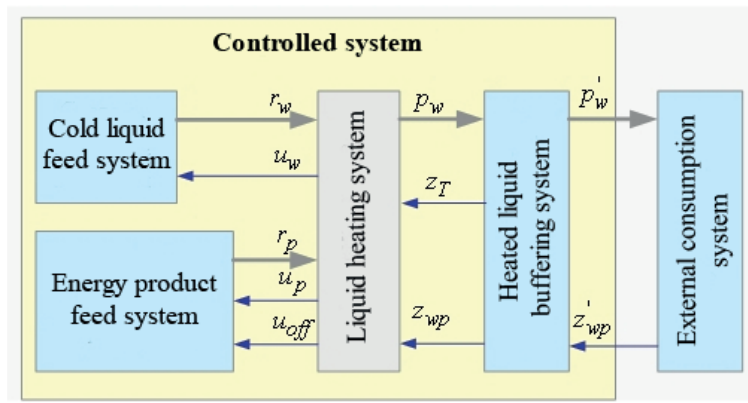

Fig. 1. Controlled liquid heating system with the opportunity of target operation optimization

Verification of the model performance was implemented in the software constructor EFFLI. Necessary mechanisms and systems needed for the synthesized conversion system functioning were developed for each element of Petri net.

## 5. The interface model of the system in the form of EFFLI objects

Representation of the model of the controlled system in the form of EFFLI objects with named inputs and outputs allows to build an interface model (scheme) of the controlled system (Fig. 3).

The object of the EFFLI environment (system, subsystem, mechanism) has an object-oriented structure.

Object-oriented approach lies in the fact that each object "knows" the way to perform its specialized function using specialized mechanism. Specialized mechanism extracts all necessary information products from the input port sections, processes them and sends to an output section or output sections of the port.

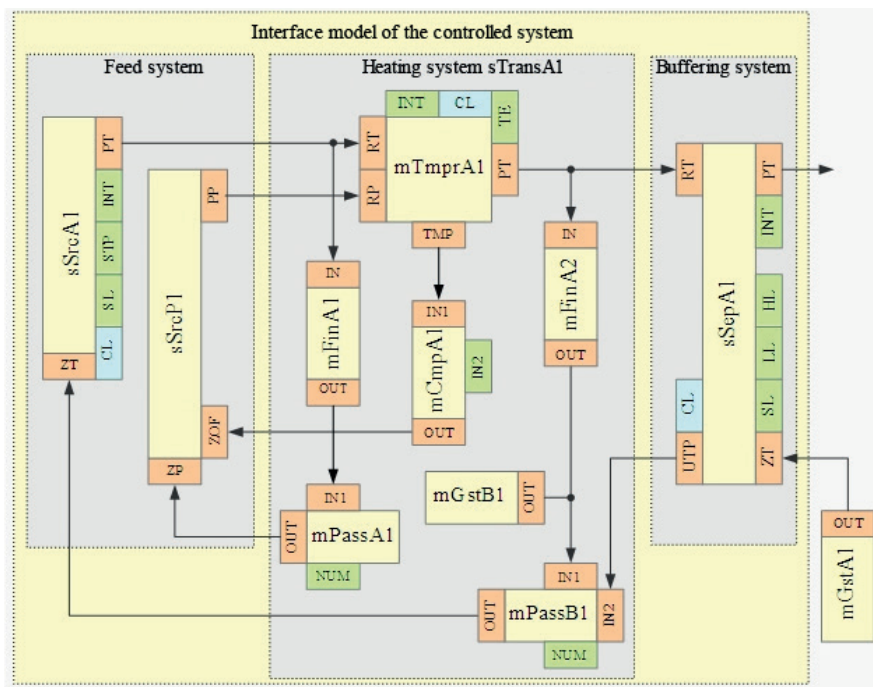

Fig. 3. The interface model of the controlled heating system

Thus, each arrow of the interface model is, in fact, the object that provides the functions of reception, movement and release of information product. Each section of the port performs the functions of reception, buffering and release of the information product.

The order of service of EFFLI objects by the operating system is specified in the program settings.

### 6. Description of EFFLI objects

#### 6. 1. The sSepA buffering system

The sSepA buffering system is the object of the controlled heating system and provides the opportunity of independent functioning of the batch heating system (Fig. 4). Figure 1 at the end of the name of the sSepA1 system in the connection diagram means that this is the first instance of such object in the overall model diagram.

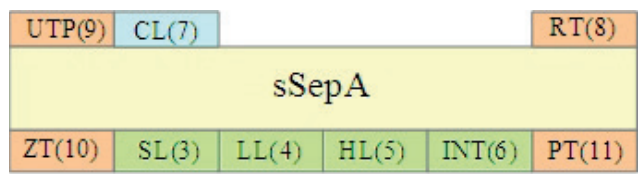

Fig. 4. The interface model of the buffering system in the form of the EFFLI object

Assignment of the sSepA1 port sections is given in Table 1.
Each section of the port has its own address that is listed in brackets (Fig. 4) after the port section notation.

Necessary initial parameters of the system are set before the start of the simulation in the SL, LL, HL and INT sections of the first object instance. The initial liquid stock level of internal buffering mechanism (reservoir) is set in the SL section. If the initial stock level is above a predetermined lower level in the LL section, the sSepA1 system will not generate a control signal for restocking until the current stock level falls below a predetermined level in the LL section.

Section CL is designed to monitor the current stock level of the system.

Changing any value in the instance section leads to certain changes in the operation of both the sSepA system, and controlled system. This applies to all objects discussed below. Do not be afraid of experiments. An example of a finished system with the initial settings is always available for download via the link [11].

Synthesis and detailed description of the sSepA system operation will be given later.

Table 1

Assignment of the sSepA1 port sections

| Port assignment | Not. | Inst.1 |
|---|---|---|
| Initial stock level | SL | 0 |
| Lower stock level | LL | 1,2 |
| Upper stock level | HL | 2 |
| Technological product release intensity | INT | 0,01 |
| Current stock level | CL | 0 |
| Technological product feed | RT | 0 |
| Restocking request | UTP | 0 |
| Assignment for release of the required volume of technological product | ZT | 0 |
| Technological product release | PT | 0 |

#### 6. 2. The sSrcA buffering system

The SSrcA buffering system is a modification of the sSepA system with limited functionality (Fig. 5). There is no possibility for technological product feed. This is useful when there is no need to explicitly specify the initial source of the technological product. The required number of operation cycles is specified by the initial stock level.

Assignment of the sSrcA1 port sections is given in Table 2.

A distinctive feature of the sSrcA1 system is the availability of the STP section. In the presence of "1" in the port section, the simulation process stops at the moment of achieving the zero stock level.

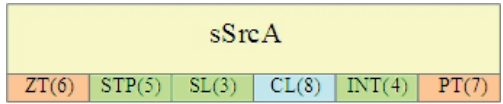

Fig. 5. The interface model of the buffering system with limited functionality by the input in the form of the EFFLI object

Table 2

Assignment of the sSrcA1 port sections

| Port assignment | Not. | Inst.1 |
|---|---|---|
| Initial stock level | SL | 30 |
| Technological product release intensity | INT | 0,01 |
| Stop of the system when achieving the zero stock level | STP | 0 |
| Assignment for the release of the required volume of technological product | ZT | 0 |
| Technological product release | PT | 0 |
| Current stock level | CL | 0 |

### 6. 3. The sSrcP energy product feed system

The model of the sSrcP energy product feed system is a model with unlimited stock (Fig. 6).

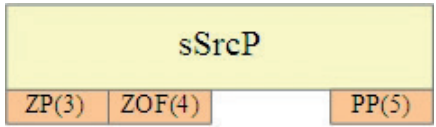

Fig. 6. The interface model of the energy product feed system in the form of the EFFLI object

Feed of the impulse signal with a certain amplitude to the ZP section provides a start of the energy product feed to the PP section with a quantity, equal to the impulse signal amplitude.

Feed of the unit impulse signal to the ZOF section provides a stop of the energy product feed.

Assignment of the sSrcP system port sections is given in Table 3.

Table 3

Assignment of the sSrcP system port sections

| Port assignment | Not. | Inst.1 |
|---|---|---|
| Assignment for the start of the energy product feed | ZP | 0 |
| Signal of the stop of the energy product feed | ZOF | 0 |
| Technological product release | PP | 0 |

### 6. 4. Conversion system mechanisms

The conversion system consists of objects, each also being the system. But these are systems, which do not have the degrees of freedom. Therefore, they are defined in the paper as mechanisms. Let us consider the operation principles of these mechanisms.

### The mTmprA mechanism

The mTmprA heating mechanism provides the heating process simulation based on the initial ambient temperature, heating tank geometry and material. Interface model of the heating mechanism is shown in Fig. 7.

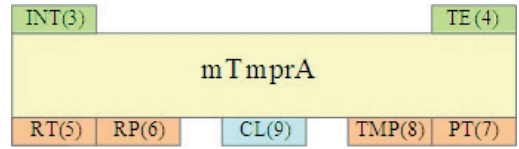

Fig. 7. The interface model of the heating mechanism in the form of the EFFLI object

Assignment of the mTmprA mechanism port sections is shown in Table 4.

Table 4

Assignment of the mTmprA mechanism port sections

| Port assignment | Not. | Inst.1 |
|---|---|---|
| Technological product release intensity | INT | 0,01 |
| Ambient temperature | TE | 20 |
| Technological product feed | RT | 0 |
| Energy product feed | RP | 0 |
| Technological product release | PT | 0 |
| Current heating temperature | TMP | 0 |
| Tank loading level | CL | 0 |

Model of the heating-cooling process is well known and tested in practice [12]. The current temperature of the tank with liquid is changed in accordance with the expression

$$T_c = T_{c-1} + dT_g - dT_d,$$

where $T_c$ – current temperature value; $T_{c-1}$ – previous temperature value; $dT_g$ – temperature growth; $dT_d$ – temperature decrease.

In turn,

$$dT_g = \frac{P \cdot dt + T_E}{c_v m_v + c_w m_w}; \quad dT_d = \frac{\eta \cdot s(T_c - T_E)dt}{c_v m_v d_v},$$

where $P$ – energy source power; $dt$ – temperature change; $T_E$ – ambient temperature; $c_v$ – tank heat capacity; $c_w$ – liquid heat capacity; $m_v$ – tank mass; $m_w$ – liquid mass; $\eta$ – ambient heat conduction; $s$ – area of the tank walls; $d_v$ – thickness of the tank walls.

Determining the liquid mass and setting the dimensional parameters, which define the tank mass, and also setting the radiation source power, it is possible to get a clear picture of the heating process physics depending on various factors.

Timing diagrams of the mTmprA mechanism are conceptually shown in Fig. 8, assuming that the process is shown in a continuous coordinate system. Hereinafter, the notation of the x-axis corresponds to the name of the port section.

Description of the object functioning is made based on the assumption that physical objects and information products move within the system.

Portion of cold liquid is supplied to the RT port section. Section CL allows to control the current liquid level in the buffer tank. Since the energy product feed through the RT section, liquid heating starts and stops when achieving a predetermined heating temperature in the TE section by stopping the energy product feed. After heating, the liquid is transferred to the heating mechanism output through the PT section. Information about the current heating temperature enters the TMP section.

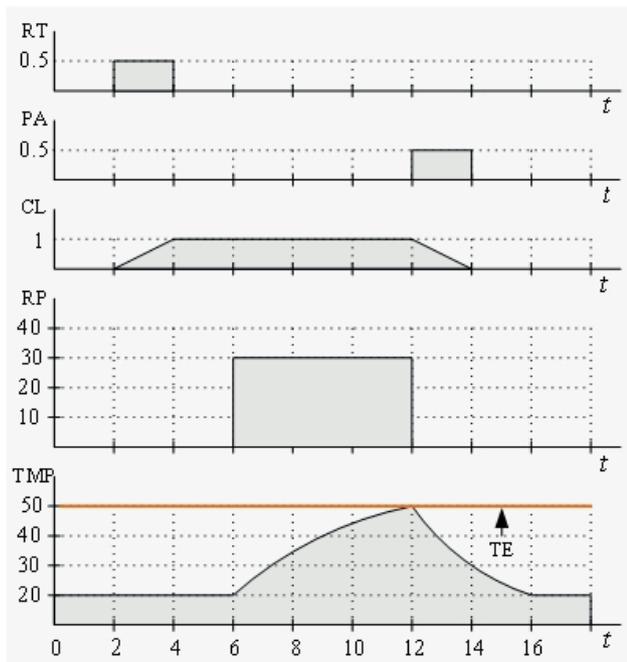

Fig. 8. Timing diagrams of the heating mechanism operation

**The mPassB mechanism**

The mPassB mechanism ensures a coordinated transfer of the NUM section data, subject to the prior receipt of unit signals to the IN1 and IN2 section. The interface model of the mPassB mechanism is shown in Fig. 9.

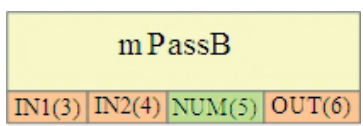

Fig. 9. The interface model of the heating mechanism in the form of the EFFLI object

Assignment of the mPassB mechanism port sections is shown in Table 5.

Table 5

Assignment of the mPassB mechanism port sections

| Port assignment | Not. | Inst.1 |
|---|---|---|
| Unit-level input signal 1 | IN1 | 0 |
| Unit-level input signal 2 | IN2 | 0 |
| Input signal 3 | NUM | 1 |
| Output signal | OUT | 0 |

Fig. 10 shows the timing diagrams of the mPassB mechanism operation.

The information to be transmitted in the OUT section before startup of the controlled system is written in the Num section. Unit impulse signals that are fed to the IN1 and IN2 sections, are transferred to the mem1 and mem2 internal memory cells. Once the values of both internal memory cells become non-zero, the value in the NUM section is transmitted to the OUT port section. Herewith, memory cells are reset.

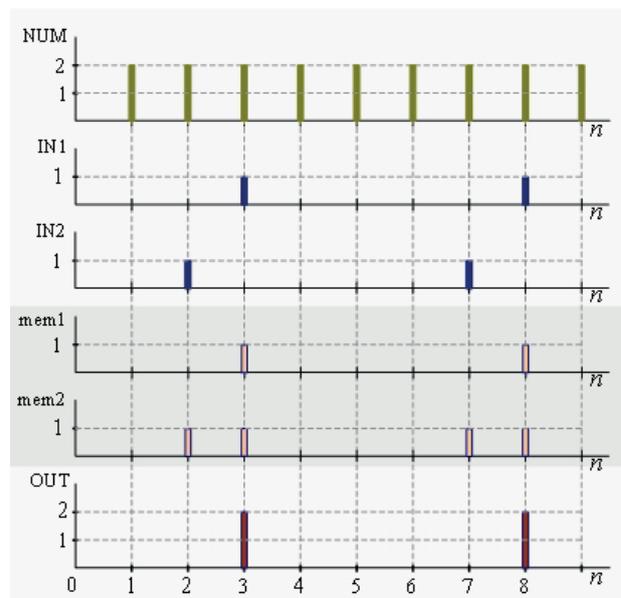

Fig. 10. The timing diagrams of the mPassB mechanism operation

The mPassA mechanism operation is not fundamentally different from the mPassB. It has only one coordinating signal, which is fed to the IN section and one mem memory cell.

**The mFinA mechanism**

The mFinA mechanism determines the time of the completion of technological product transfer. The interface model of the mFinA mechanism is shown in Fig. 11.

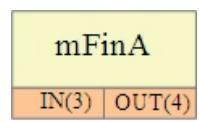

Fig. 11. The interface model of the mFinA mechanism in the form of the EFFLI object

Assignment of the mFinA mechanism port sections is shown in Table 6.

Table 6

Assignment of the mFinA mechanism port sections

| Port assignment | Not. | Inst.1 |
|---|---|---|
| Input signal | IN | 0 |
| Output signal | OUT | 0 |

Control of the time of termination of liquid feed or release is ensured by the mFinA mechanism. The principle of its operation is based on the fact that the in the presence of non-zero flow in the unit time interval, a high logic level signal is recorded in the mem memory cell. Termination of the liquid flow leads to a situation where the unit-level memory cell corresponds to zero liquid flow. This combination leads to the cleaning of the memory cell and the formation of a high-level impulse signal and its transmission to the OUT section (Fig. 12).

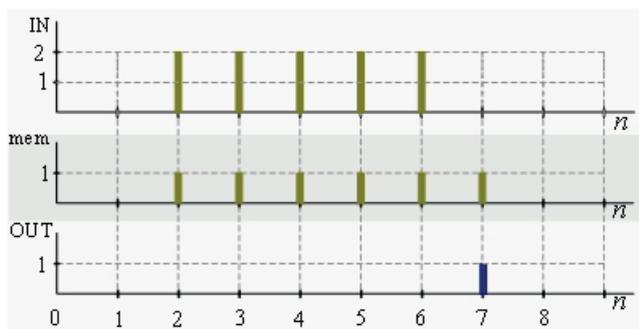

Fig. 12. The timing diagrams of the mFinA mechanism operation

**The mCmpA signal comparison mechanism (comparator)**

Comparison of two signals is provided by the mCmpA mechanism. Its interface model is shown in Fig. 13.

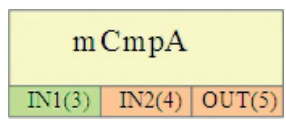

Fig. 13. The interface model of the mCmpA signal comparison mechanism in the form of the EFFLI object

Assignment of the mCmpA mechanism port sections is shown in Table 7.

Table 7

Assignment of the mCmpA mechanism port sections

| Port assignment | Not. | Inst.1 |
|---|---|---|
| Input signal | IN1 | 50 |
| Input signal | IN2 | 0 |
| Output signal | OUT | 0 |

The mCmpA comparator is a device for comparing signals in IN1 and IN2 sections. At a time of equality of these signals, unit impulse signal is transmitted to the OUT section (Fig. 14).

The mGstB mechanism ensures the formation of unit impulse signal at the time of the controlled system startup, and mGstA mechanism emulates the feed of the signal of external consumption system demand in a certain amount of heated liquid. In this case, the mechanism generates impulse of predetermined magnitude per unit of the sampling interval. The appropriate volume of the heated liquid is fed to the output through the PT section of the sSepA system.

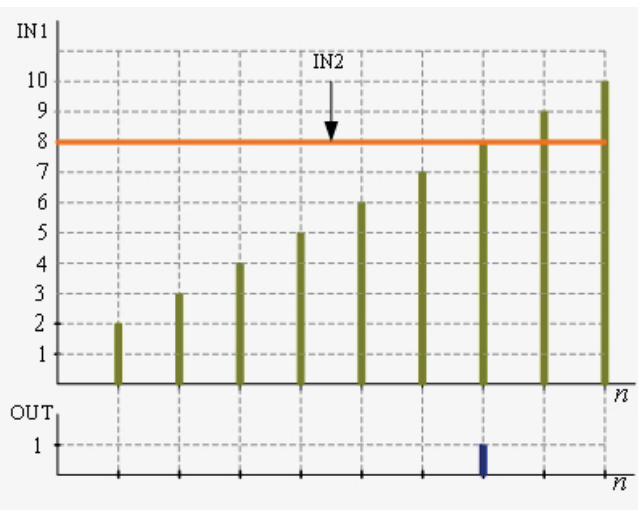

Fig. 14. The timing diagrams of the mCmpA mechanism operation

## 7. The results of the controlled system startup in the EFFLI environment

The EFFLI software constructor is developed in the EXCEL application using the programming environment Visual Basic for Applications. It's free, so the model is available to any user of Microsoft Office standard package. It is necessary first to set the average safety level and not to disable macros at the program startup.

Ready-to-startup model of the controlled system is available on the resource [11]. It is possible to monitor the process of the system operation by setting different energy product feed rate parameters, and in the pages of reporting, build or see a registration model of one or another operation depending on the current task.

The result of startup of the assembled system is shown in Fig. 15. There are timing diagrams, associated with the change in temperature (Fig. 15, *a*) and technological product stock in the mTmprA1 mechanism (Fig. 15, *b*), and Fig. 16 shows timing diagrams, associated with the heated liquid level in the tank (Fig. 16, *a*) and portions of technological product restocking in the sSepA1 system (Fig. 15, *b*).

Independent change in the volume of the liquid portion, fed to the conversion system input, the change in the energy product feed rate, change in the liquid safety stock level in the buffering system and the upper boundary of the liquid level is possible within the controlled system.

Thus, this structure of the controlled system has a complete set of features for implementing a full-scale parametric optimization based on this structure.

Studies have shown that the developed structure of the controlled system allows to implement an operation mode with continuous feed of cold liquid, its heating and release by selecting the energy product rate per time unit. Thus, as expected, heating process optimization possibilities are lost completely. A certain energy product feed rate that provides the specified temperature of the heated liquid corresponds to each feed – release rate.

This means that rigid systems with continuous feed – release of raw product are a special case of completely-controllable systems with architecture that provides the optimal control possibility.

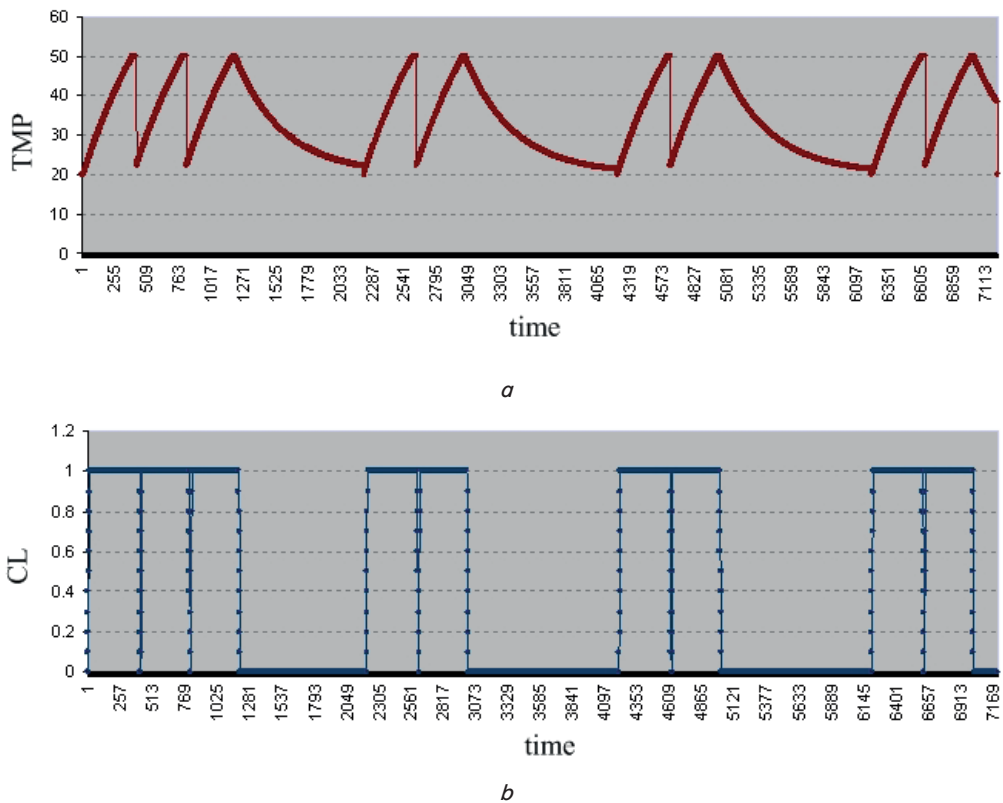

Fig. 15. Timing diagrams of technological process of heating mechanism: *a* — change in the tank temperature; *b* — liquid level in the tank of the mTmprA1heating mechanism

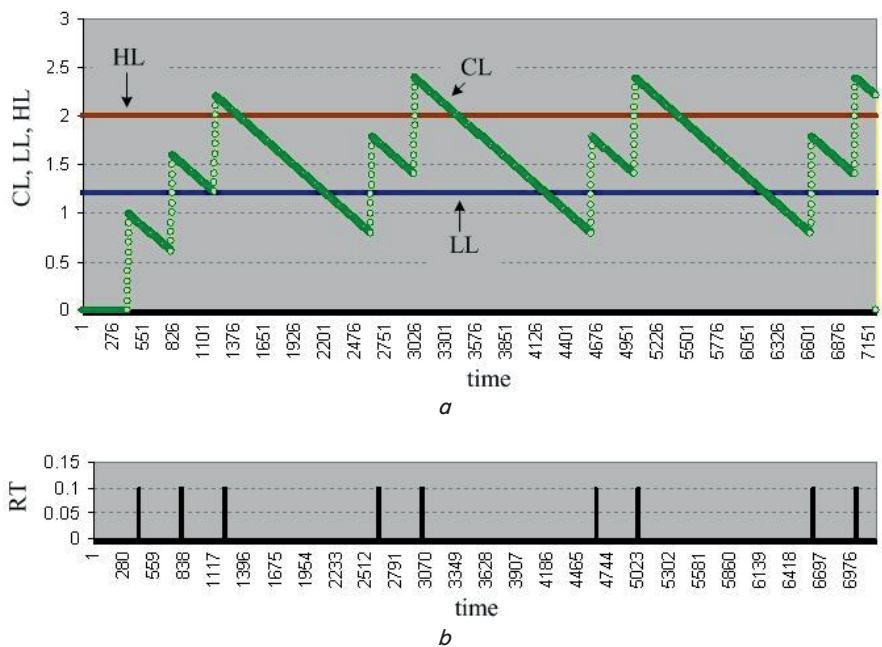

Fig. 16. Timing diagrams of technological process of buffering system: *a* — change in the heated liquid stock level; *b* — heated liquid restocking in the sSepA1 buffering system

## 8. Conclusions

The internal structure of the conversion system with batch feed of one technological product is synthesized. The architecture of the controlled system with the maximum possible number of degrees of freedom is developed.

It is found that the internal structure of the conversion system with batch feed – release of raw products includes:
– raw, energy and information product movement mechanisms
– technological raw product conversion mechanism
– raw product receipt completion registration mechanism,

– technological process quality control mechanism,

– conversion process completion control mechanism,

– finished product release completion registration mechanism,

– raw product control signal feed synchronization mechanism,

– energy product control signal feed synchronization mechanism.

It is found that the conversion system with batch feed – release of raw products in combination with technological product buffering system provides the possibility of an independent energy product feed change during the independent installation in the buffering system of safety stock and the upper stock level. This allows to transfer a finished product with desired consumer properties and in the required amount to the consumption system.

A feature of the controlled system architecture is that the functions of ensuring the specified quality of the consumer product and features of ensuring the release of the consumer product in the required volume are divided between specialized systems (conversion system and buffering system). It is this structure that allows potential implementation of optimal control.

It is experimentally found that systems with continuous feed – release of raw product are a special case of fully controllable systems with the architecture that provides the optimal control possibility.

## References

1. Gavrilov, D. A. Upravlenie proizvodstvom na baze standarta MRP [Text] / D. A. Gavrilov. – Izdatelskiy dom «Piter», 2002. – 320 p.
2. Kirk, E. Optimal Control Theory: An Introduction (Dover Books on Electrical Engineering) [Text] / E. Kirk. – Dover Publications. – 2004. – 464 p.
3. Athans, M. Optimal Control: An Introduction to the Theory and Its Applications [Text] / M. Athans, L. Falb. – Dover Publications, 2006. – 879 p.
4. Zhang, S. Optimal Control Strategy Design Based on Dynamic Programming for a Dual-Motor Coupling-Propulsion System [Text] / S. Zhang, C. Zhang, G. Han, Q. Wang // The Scientific World Journal. – 2014. – Vol. 2014. – P. 1–9. doi: 10.1155/2014/958239
5. Pierre, A. Optimization Theory with Applications [Text] / A. Pierre. – Donald Courier Dover Publications. 1986. – 612 p.
6. Everett, E. Operations Change Interactions in a Service Environment: Attitudes, Behaviors, and Profitability [Text] / E. Everett // Journal of Operations Management. – 1981. – Vol. 2, Issue 1. – P.63-76. doi: 10.1016/0272-6963(81)90036-x
7. Frederic, D. Systems modeling: analysis and operations research [Text] / D. Frederic // Modeling and Simulation Fundamentals: Published Online. – 2010. –Vol. 6. – P. 147–180. doi: 10.1002/9780470590621.ch6
8. Crassidis, L. Optimal Estimation of Dynamic Systems [Text] / L. J. Crassidis, L. J. Junkins. – E-book Google, 2004. – 608 p.
9. Kagramanyan, S. L. Modelirovanie i upravlenie gornorudnyimi predpriyatiyami [Text] / S. L. Kagramanyan, A. S. Davidkovich, V. A. Malyishev, O. Burenzhargal, Sh. Otgonbileg. – Nedra, 1989. – 360 p.
10. Lutsenko, I. A. A practical approach to selecting optimal control criteria [Text] / I. A. Lutsenko // Technology audit and production reserves. – 2014. – Vol. 2/1(16). – P. 32–35. Available at: http://journals.uran.ua/tarp/article/view/23432/20906
11. Lutsenko, I. A. Samples [Electronic resource] / I. A. Lutsenko. – Krivoy Rog, 2014. – Available at: http://uk.effli.info/index.php/systems-engineering-samples
12. Amelkin, V. V. Differentsialnyie uravneniya v prilozheniyah [Text] / V. V. Amelkin. – Nauka, 1987. – 160 p.